\documentclass[12pt]{article} 
\usepackage[utf8]{inputenc} 
\usepackage{geometry}
\usepackage{graphicx} 
\usepackage{amssymb}
\usepackage{amsmath}
\usepackage{cite}
\usepackage{mathrsfs}
\usepackage{booktabs} 
\usepackage{array} 
\usepackage{paralist}
\usepackage{verbatim} 
\usepackage{subfig}
\usepackage{fancyhdr} 
\usepackage{sectsty}
\allsectionsfont{\sffamily\mdseries\upshape} 
\usepackage[nottoc,notlof,notlot]{tocbibind} 
\usepackage[titles,subfigure]{tocloft} 

\title{Exact solutions of the boundary-value problems for the Helmholtz equation in a layer with polynomials in the right-hand sides of the equation and of the boundary conditions}
\author{\bf Oleg D. Algazin}
\date{Bauman Moscow State Technical University, Moscow, Russia} 

\begin{document}
\maketitle
\thispagestyle{empty}

\begin{abstract}
         In a multidimensional infinite layer bounded by two hyperplanes, the inhomogeneous Helmholtz equation with a polynomial right-hand side is considered. It is shown that the Dirichlet and Dirichlet-Neumann boundary-value problems with polynomials in the right-hand sides of the boundary conditions have a solution that is a quasipolynomial that contains, in addition to the power functions, also hyperbolic or trigonometric functions. This solution is unique in the class of slow growth functions if the parameter of the equation is not an eigenvalue. An algorithm for constructing this solution is given and examples are considered.
\end{abstract}


\section*{Introduction}

 Many problems of mathematical physics lead to the Helmholtz equation, for example, problems associated with steady-state oscillations (mechanical, acoustic, electromagnetic, etc.), and the problems of  diffusion of some gases in the presence of decay or chain reactions. Also, any equation of elliptic type with constant coefficients is reduced to the Helmholtz equation \cite{TS}.

     In this paper, exact solutions are obtained in the form of quasipolynomials of the Dirichlet and Dirichlet-Neumann boundary value problems for the Helmholtz equation in a layer in the case when the right-hand side of the Helmholtz equation and the right-hand sides of the boundary conditions are polynomials. If the parameter of the Helmholtz equation tends to zero, then the Helmholtz equations turn into the Poisson equation, and the quasipolynomial solutions of the boundary value problems turn into polynomial solutions of the boundary value problems for the Poisson equation \cite{Alg1}. In the same way, exact polynomial solutions of  boundary value problems for the Tricomi equation in a strip are obtained\cite{Alg2},\cite{Alg3}. The search for solutions of partial differential equations in the form of polynomials or quasipolynomials has been the subject of work by many authors \cite{Nik} -\cite{Nyt}.
      
\section{Statement of the problem.  }

We will  consider the inhomogeneous Helmholtz equation  with a polynomial right-hand side in an unbounded domain (in the layer) 
\[
                      \Delta u(x,y)+\nu u(x,y)=P(x,y),~~\nu \in \mathbb R,~~x\in \mathbb R^n,~~0<y<a,   \eqno{(1)}
\]
 where $x=(x_1,…,x_n ),~\Delta$ is the Laplace operator,
\[
\Delta=\frac{\partial^2 }{\partial x_1^2}+...+\frac{\partial^2 }{\partial x_n^2}+\frac{\partial^2 }{\partial y^2},
\]
$P(x,y)$ is a  polynomial in the variables $x$ and $y$.

On the boundary of the layer we set the Dirichlet boundary conditions
\[
u(x,0)=\phi(x),~~u(x,a)=\psi(x),~~x\in\mathbb{R}^n, \eqno{(2)}
\]
where $\phi(x),\psi(x)$ are polynomials.

     Further in section 2 it is shown that the inhomogeneous Helmholtz equation (1) has polynomial solutions and the formula for obtaining such a solution is given. If $\tilde u(x,y)$ is some polynomial solution of the inhomogeneous Helmholtz equation (1), then for the function $v(x,y)=u(x,y)-\tilde u(x,y)$ we obtain a homogeneous Helmholtz equation
\[
                      \Delta u(x,y)+\nu u(x,y)=0,~~\nu \in \mathbb R,~~x\in \mathbb R^n,~~0<y<a,   \eqno{(3)}
\]
            and Dirichlet boundary conditions
\[
v(x,0)=\phi(x)-\tilde u(x,0),~~v(x,a)=\psi(x)-\tilde u(x,a),~~x\in\mathbb{R}^n.\eqno{(4)}
\]
     If we solve the Dirichlet problem for the homogeneous Helmholtz equation (3), (4), we obtain the solution of the Dirichlet problem for the inhomogeneous Helmholtz equation (1), (2) by the formula
\[
u(x,y)=v(x,y)+\tilde u(x,y).
\]
     The mixed Dirichlet-Neumann boundary value problem with boundary conditions 
\[
u(x,0)=\phi(x),~~u_y(x,a)=\psi(x),~~x\in\mathbb{R}^n, \eqno{(5)}
\]
is considered in a similar way.
This problem also reduces to the problem for homogeneous equation (3) with boundary conditions
\[
v(x,0)=\phi(x)-\tilde u(x,0),~~v_y(x,a)=\psi(x)-\tilde u_y(x,a),~~x\in\mathbb{R}^n.
\]

      Solutions $u(x,y)$ of the boundary value problems (1), (2) and (1), (5) we will  sought in the class of functions of slow growth with respect to the variable $x$ for each fixed $y$ from the interval $(0,a)$, i.e. for $\forall y \in(0,a)$ there exists $m\ge0$ such that
\[
\int_{\mathbb{R}^n}|u(x,y)|(1+|x|^2)^{-m}dx<\infty,~~|x|=\sqrt{x_1^2+...+x_n^2}.\eqno{(6)}
\]
                                                 Therefore, one can apply the Fourier transform for generalized functions of slow growth with respect to the variable $x$ \cite{Vla}.

\section{The polynomial solution of the inhomogeneous Helmholtz equation}

     The inhomogeneous Helmholtz equation (1) with the polynomial right-hand side $P(x,y)$,
\[
                      \Delta u(x,y)+\nu u(x,y)=P(x,y),~~\nu \in \mathbb R,~~x\in \mathbb R^n,~~y\in \mathbb R,
\]
has polynomial solutions, one of which can be obtained by the following directly verified formula.
     For $\nu \ne0$
\[
                              u(x,y)=\frac{P(x,y)}{\nu}+\sum_{j=1}^{[k/2]}\frac{(-1)^j}{\nu^{j+1}}\Delta^j P(x,y)  , \eqno{(7)}                                                  \] 
where $k$  is the largest of the degrees of the monomials of the polynomial $P(x,y),  [k/2]$ is the integer part of the number $k/2$.   
      For $\nu=0$, we have the Poisson equation and its polynomial solutions are given in \cite{Alg1}.

\textbf{Example 1. }
\[
x=(x_1,x_2 ),P(x,y)=3x_1^2 x_2 y^2+5x_1 x_2^2 y,~~k=5,~~[k/2]=2.
\]
By the formula (7) we obtain
\[
u(x,y)=\frac{1}{\nu} (3x_1^2 x_2 y^2+5x_1 x_2^2 y)-\frac{1}{\nu^2}(6x_2 y^2+10x_1 y+6x_1^2 x_2 )+\frac{1}{\nu^3}24x_2.
\]

\section{	Solution of the Dirichlet problem for the Helmholtz equation in the case $\nu=-\lambda^2$}

       Since the solution of the Dirichlet problem for an inhomogeneous equation reduces to the solution of the Dirichlet problem for a homogeneous equation, we will consider the Dirichlet problem for a homogeneous equation
\[
\Delta u(x,y)-\lambda^2 u(x,y)=0,~~\lambda>0,~~x\in \mathbb R^n,~~0<y<a,   \eqno{(8)}
\]
\[
u(x,0)=\phi(x),~~u(x,a)=\psi(x),~~x\in\mathbb{R}^n, \eqno{(9)}
\]
where $\phi(x),\psi(x)$ are polynomials.

We will apply the Fourier transform with respect to $x$ [10]:
\[
\mathscr{F}_x\left[u(x,y)\right](t,y)=U(t,y),~~\mathscr{F}_x\left[\phi(x)\right](t)=\Phi(t),~~\mathscr{F}_x\left[\psi(x)\right](t)=\Psi(t)
\]
and we obtain the boundary value problem for an ordinary differential equation of the second order with parameter $t\in\mathbb{R}^n$,
\[
-(\lambda^2+|t|^2)U(t,y)+U_{yy} (t,y)=0,~~t\in \mathbb R^n,~~0<y<a,  \eqno{(10)}                                 
\]
\[
          U(t,0)=\Phi(t),~~U(t,a)=\Psi(t)  .                                                            \eqno{(11)}
\]
The unique solution to the boundary value problem (10), (11) is given by the formula
\[
U(t,y)=L_n (|t|,a-y)\Phi(t)+L_n (|t|,y)\Psi(t),                            \eqno{(12)}
\]
where
\[
L_n(|t|,y)=\frac{\sinh(y\sqrt{|t|^2+\lambda^2})}{\sinh(a\sqrt{|t|^2+\lambda^2})}.
\]
Applying the inverse Fourier transform, we obtain the unique solution of the Dirichlet problem (8), (9) in the class of functions of slow growth in the form of convolution
\[
                                       u(x,y)=l_n (|x|,a-y)*\phi(x)+l_n (|x|,y)*\psi(x),                           \eqno{(13)}
\]
where
\[
  l_n (|x|,y)=\mathscr{F}_t^{-1}\left[L_n (|t|,y)\right](x,y).
\]
To find the convolution (13) with polynomials $\phi(x)$ and $\psi(x)$ it suffices to consider the case of a monomial.

	\subsection{Case $n = 1$}
We first consider the flat case, $n=1, ~  x\in \mathbb R, ~  (x,y)\in\mathbb R \times(0,a)$.

$L_1 (|t|,y)=L_1 (t,y)$ is even  function of the variable $t$,
\[ 
L_1(|t|,y)=L_1(t,y)=\frac{\sinh(y\sqrt{t^2+\lambda^2})}{\sinh(a\sqrt{t^2+\lambda^2})}.
\]
Let $\phi(x)=0,~~\psi(x)=x^0=1$. 
Solution to the Dirichlet problem is
\begin{gather*}
u_0(x,y)=l_1(x,y)*\psi(x)=\int_{-\infty}^{\infty}\psi(x-t)l_1(t,y)dt=\int_{-\infty}^{\infty}l_1(t,y)dt=
\lim_{x\to 0}\int_{-\infty}^{\infty}l_1(t,y)e^{ixt}dt=\\
=\lim_{x\to 0}\mathscr{F}_t\left[l_1(t,y)\right](x,y)=\lim_{x\to 0}L_1(x,y)=\lim_{x\to 0}\frac{\sinh(y\sqrt{x^2+\lambda^2})}{\sinh(a\sqrt{x^2+\lambda^2})}=\frac{\sinh(\lambda y)}{\sinh(\lambda a)}.
\end{gather*}

Now let $\phi(x)=0,~\psi(x)=x^k$.
The corresponding solution to the Dirichlet problem is
\begin{gather*}
u_k(x,y)=l_1(x,y)*\psi(x)=\int_{-\infty}^{\infty}\psi(x-t)l_1(t,y)dt=\int_{-\infty}^{\infty}(x-t)^k l_1(t,y)dt=\\
=\int_{-\infty}^{\infty}\sum_{j=0}^kC_k^jx^{k-j}t^j(-1)^jl_1(t,y)dt=\sum_{j=0}^kC_k^jx^{k-j}(-1)^j\int_{-\infty}^{\infty}t^jl_1(t,y)dt,
\end{gather*}
where $C_k^j=k!/j!(k-j)!$  are binomial coefficients. Since the last integral for odd $j$ is equal to zero due to the parity of $l_1 (t,y)$ with respect to the variable $t$, then
\[
u_k(x,y)=\sum_{m=0}^{[k/2]}C_k^{2m}x^{k-2m}\int_{-\infty}^{\infty}t^{2m}l_1(t,y)dt=\sum_{m=0}^{[k/2]}C_k^{2m}x^{k-2m}p_{2m}(\lambda y),
\]
where $[k/2]$ is the integer part of the number $k/2$.  Using the properties of the Fourier transform, we obtain
 \[
p_{2m}(\lambda y)=\int_{-\infty}^{\infty}t^{2m}l_1(t,y)dt=\lim_{x\to 0}\int_{-\infty}^{\infty}t^{2m}l_1(t,y)e^{ixt}dt=
\]
\[
=\lim_{x\to 0}\mathscr{F}_t\left[t^{2m}l_1(t,y)\right](x,y)=(-1)^m\lim_{x\to 0}\frac{\partial^{2m}}{\partial x^{2m}}L_1(x,y).
\]
The functions $p_{2m} (\lambda,y)$ are coefficients of expansion in a power series in $x$ of the function
\[
L_1(x,y)=\frac{\sinh(y\sqrt{x^2+\lambda^2})}{\sinh(a\sqrt{x^2+\lambda^2})}=\sum_{m=0}^{\infty}p_{2m}(\lambda,y)\frac{(-1)^mx^{2m}}{(2m)!},
\]
i.e. $L_1 (x,y)$ is the generating function for $p_{2m} (\lambda,y)$. These functions can be calculated using the recurrent formula
\[
p_0(\lambda,y)=\frac{\sin(\lambda y)}{\sin(\lambda a)},~~p_{2m}(\lambda,y)=-(2m-1)\frac{1}{\lambda}\frac{\partial}{\partial\lambda}p_{2m-2}(\lambda,y),~~m=1,2,\dots\eqno{(14)}
\]
Let us prove this formula.

       Because the
\[
f(s)=\frac{\sinh(ys)}{\sinh(ya)},~~s=\sqrt{x^2+\lambda^2},
\]
is an even analytic function of the complex variable $s$, and its singular points $\pm i\pi k/a,~k\in\mathbb N$ lie on the imaginary axis, then in the circle $|s|<\pi/a$ we have the expansion
\[
f(s)=\sum_{n=0}^{\infty}a_{2n}s^{2n}=\sum_{n=0}^{\infty}a_{2n}(x^2+\lambda^2)^n=\sum_{n=0}^{\infty}a_{2n}\sum_{m=0}^nC_n^m\lambda^{2n-2m}x^{2m}=
\]
\[
=\sum_{m=0}^{\infty}x^{2m}\sum_{n=m}^{\infty}a_{2n}C_n^m\lambda^{2n-2m}.
\]
Hence,
\[
p_{2m}(\lambda,y)=(-1)^m(2m)!\sum_{n=m}^{\infty}a_{2n}C_n^m\lambda^{2n-2m},
\]
\[
p_0(\lambda,y)=\sum_{n=0}^{\infty}a_{2n}\lambda^{2n}=f(\lambda)=\frac{\sinh(\lambda y)}{\sinh(\lambda a)},
\]
\[
p_{2m-2}(\lambda,y)=(-1)^{m-1}(2m-2)!\sum_{n=m-1}^{\infty}a_{2n}C_n^{m-1}\lambda^{2n-2m+2},
\]
\[
-(2m-1)\frac{1}{\lambda}\frac{\partial}{\partial \lambda}   p_{2m-2} (\lambda,y))=(-1)^m(2m)!\sum_{n=m}^{\infty}a_{2n}\frac{2n-2m+2}{2m}C_n^{m-1}\lambda^{2n-2m}=
\]
\[
=(-1)^m(2m)!\sum_{n=m}^{\infty}a_{2n}C_n^m\lambda^{2n-2m}=p_{2m}(\lambda,y),
\]
q.e.d.

Thus,
\[
p_{2m}(\lambda,y)=(2m-1)!!\left(-\frac{1}{\lambda}\frac{\partial}{\partial\lambda}\right)^m\frac{\sinh(\lambda y)}{\sinh(\lambda a)}.
\]
For example,
\begin{gather*}
p_2(\lambda,y)=-\frac{y\cosh(\lambda y)}{\lambda \sinh(\lambda a)}+\frac{a\sinh(\lambda y)\cosh(\lambda a)}{\lambda\sinh^2(\lambda a)},\\
p_4(\lambda,y)=-\frac{3y\cosh(\lambda y)}{\lambda^3 \sinh(\lambda a)}+\frac{3y^2\sinh(\lambda y)}{\lambda^2 \sinh(\lambda a)}-\frac{6ya\cosh(\lambda y)\cosh(\lambda a)}{\lambda^2 \sinh^2(\lambda a)}+\frac{3a\sinh(\lambda y)\cosh(\lambda a)}{\lambda^3 \sinh^2(\lambda a)}+\\+\frac{6a^2\sinh(\lambda y)\cosh(\lambda a)^2}{\lambda^2 \sinh^3(\lambda a)}-\frac{3a^2\sinh(\lambda y)}{\lambda^2\sinh(\lambda a)}.
\end{gather*}

As $\lambda$ tends to zero, the functions $p_{2m}(\lambda,y)$ go over into polynomials $p_{2m}(y)$, which were considered in \cite{Alg1}. For example,
\[
\lim_{\lambda\to 0}p_0(\lambda y)=\lim_{\lambda\to 0}\frac{\sinh(\lambda y)}{\sinh(\lambda a)}=\frac{y}{a},
\]
\[
\lim_{\lambda\to 0}p_2(\lambda y)=-\frac{y}{3a}(y^2-a^2),~~\lim_{\lambda\to 0}p_4(\lambda y)=\frac{y}{15a}(3y^4-10y^2a^2+7a^4).
\]

We write down the first few solutions of the Dirichlet problem for the Helmholtz equation.
\[
u_k(x,y)=\sum_{m=0}^{[k/2]}C_k^{2m}x^{k-2m}p_{2m}(\lambda y).
\]
\[
u_0(x,y)=\frac{\sinh(\lambda y)}{\sinh(\lambda a)},~~u_1(x,y)=x\frac{\sinh(\lambda y)}{\sinh(\lambda a)},
\]
\[
u_2(x,y)=x^2\frac{\sinh(\lambda y)}{\sinh(\lambda a)}-\frac{y\cosh(\lambda y)}{\lambda \sinh(\lambda a)}+\frac{a\sinh(\lambda y)\cosh(\lambda a)}{\lambda\sinh^2(\lambda a)}.
\]

As $\lambda$ tends to zero, they go over into polynomial solutions of the Dirichlet problem for the Laplace equation \cite{Alg1}
\[
\lim_{\lambda\to 0}u_0(x,y)=\frac{y}{a},~\lim_{\lambda\to 0}u_1(x,y)=\frac{xy}{a},~\lim_{\lambda\to 0}u_2(x,y)=\frac{y}{3a}(3x^2-y^2+a^2).
\]

     The functions $v_k(x,y)=u_k(x,a-y)$ are solutions of the Helmholtz equation satisfying the boundary conditions
$v_k (x,0)=x^k,~~   v_k (x,a)=0,~~x\in\mathbb R$.

\textbf{Example 2. }

Consider the Dirichlet problem for the inhomogeneous Helmholtz equation
\[
\Delta u(x,y)-\lambda^2 u(x,y)=x^2y^2,~~-\infty<x<\infty,~~0<y<a,~~\lambda>0,   
\]
\[
u(x,0)=0,~~u(x,a)=0,~~-\infty<x<\infty.
\]
A particular solution to the inhomogeneous Helmholtz equation is the polynomial
\[
\tilde u(x,y)=-\frac{x^2y^2}{\lambda^2}-\frac{2y^2}{\lambda^4}-\frac{2x^2}{\lambda^4}-\frac{8}{\lambda^4}
\]
and the Dirichlet problem for the inhomogeneous Helmholtz equation reduces to the Dirichlet problem for the homogeneous Helmholtz equation for the function $v(x,y)=u(x,y)-\tilde u(x,y)$:
\[
\Delta v(x,y)-\lambda^2 v(x,y)=0,~~-\infty<x<\infty,~~0<y<a,~~\lambda>0,   
\]
\[
v(x,0)=u(x,0)-\tilde u(x,0)=\frac{2x^2}{\lambda^4}+\frac{8}{\lambda^6}
\]
\[
v(x,a)=u(x,a)-\tilde u(x,a)=\frac{x^2a^2}{\lambda^2}+\frac{2a^2}{\lambda^4}+\frac{2x^2}{\lambda^4}+\frac{8}{\lambda^6}.
\]
The solution to this problem will be the function
\[
v(x,y)=\frac{2}{\lambda^4}u_2(x,a-y)+\frac{8}{\lambda^6}u_0(x,a-y)+\left(\frac{a^2}{\lambda^2}+\frac{2}{\lambda^4}\right)u_2(x,y)+\left(\frac{2a^2}{\lambda^4}+\frac{8}{\lambda^6}\right)u_0(x,y).
\]
The solution to the original problem will be the function
\begin{gather*}
u(x,y)=\tilde u(x,y)+v(x,y)=\\
=-\frac{x^2y^2}{\lambda^2}-\frac{2y^2}{\lambda^4}-\frac{2x^2}{\lambda^4}-\frac{8}{\lambda^6}+\frac{a^2x^2\sinh(\lambda y)}{\lambda^2\sinh(\lambda a)}-\frac{a^2y\cosh(\lambda y)}{\lambda^3\sinh(\lambda a)}+\frac{a^3\sinh(\lambda y)\cosh(\lambda a)}{\lambda^3\sinh^2(\lambda a)}+\\
+\frac{\sinh(\lambda y)(2x^2+2a^2-2x^2\cosh(\lambda a)}{\lambda^4\sinh(\lambda a)}+\frac{2y\cosh(\lambda y)(\cosh(\lambda a)-1)}{\lambda^5\sinh(\lambda a)}+\\
+\frac{2a\sinh(\lambda y)\cosh(\lambda a)(1-\cosh(\lambda a))}{\lambda^5\sinh^2(\lambda a)}+\frac{8\sinh(\lambda y)(1-\cosh(\lambda a))}{\lambda^6\sinh(\lambda a)}+\frac{2x^2\cosh(\lambda y)}{\lambda^4}+\\
+\frac{2(a-y)\sinh(\lambda y)}{\lambda^5}+\frac{8\cosh(\lambda y)}{\lambda^6}.
\end{gather*}
As $\lambda\to 0$, this solution becomes a polynomial, 
\[
\lim_{\lambda\to 0}u(x,y)=\frac{1}{12}x^2y^4-\frac{1}{12}a^3x^2y-\frac{1}{180}y^6+\frac{1}{36}a^3y^3-\frac{1}{45}a^5y,
\]
which is the solution of the Dirichlet problem for the Poisson equation \cite{Alg1}
\[
\Delta u(x,y)=x^2y^2,~~-\infty<x<\infty,~~0<y<a,   
\]
\[
u(x,0)=0,~~u(x,a)=0,~~-\infty<x<\infty.
\]

\subsection{Case $n>1$}

Now we consider the case, $n>1,~~ x\in\mathbb R^n,~~ (x,y)\in\mathbb R^n\times(0,a)$.
If 
\[
\psi(x)=1,~~x\in\mathbb R^n,
\]
then the solution the Dirichlet problem with the boundary condition
\[
u(x,0)=0,~~u(x,a)=1,~~x\in\mathbb R^n,
\]
is the function
\begin{gather*}
u_0(x,y)=l_n(|x|,y)*\psi(x)=\int_{\mathbb R^n}\psi(x-t)l_n(|t|,y)dt=\int_{\mathbb R^n}l_n(|t|,y)dt=\\
=\lim_{x\to 0}\int_{\mathbb R^n}l_n(|t|,y)e^{ixt}dt=\lim_{x\to 0}\mathscr{F}_t\left[l_n(|t|,y)\right](x,y)=\lim_{x\to 0}L_n(|x|,y)=\\
=\lim_{x\to 0}\frac{\sinh(y\sqrt{|x|^2+\lambda^2})}{\sinh(a\sqrt{|x|^2+\lambda^2})}=\frac{\sinh(\lambda y)}{\sinh(\lambda a)}.
\end{gather*}

If
\[
\psi(x)=x^k,~~x\in\mathbb R^n,~~  k ~~\text{is a multiindex},
\]
then the solution the Dirichlet problem with the boundary condition
\[
u(x,0)=0,~~u(x,a)=x^k,~~x\in\mathbb R^n ,
\]
is the function
\[
u_k (x,y)=\int_{\mathbb R^n}(x-t)^kl_n(|t|,y)dt,
\]
where
\[
(x-t)^k=(x_1-t_1 )^{k_1 }\dots(x_n-t_n )^{k_n}                                                    \eqno{(15)}
\]
and this integral will be nonzero only for those monomials of polynomial (15) that contain $t_j$ in even degrees. Therefore
\[
u_k(x,y)=\sum_{m=0}^{[k/2]}C_k^{2m}x^{k-2m}\int_{\mathbb R^n}t^{2m}l_n(|t|,y)dt=
\]
\[
=\sum_{m=0}^{[k/2]}C_k^{2m}x^{k-2m}p_{2m}(\lambda,y).
\]
where
\[
C_k^{2m}=C_{k_1}^{2m_1}C_{k_2}^{2m_2}\dots C_{k_n}^{2m_n},~~[k/2]=([k_1/2],[k_2/2],\dots,[k_n/2]).
\]

Using the properties of the Fourier transform \cite{Vla}, we obtain
\[
p_{2m}(\lambda,y)=\int_{\mathbb R^n}t^{2m}l_n(|t|,y)dt=\lim_{x\to 0}\int_{\mathbb R^n}t^{2m}l_n(|t|,y)e^{ixt}dt=
\]
\[
=\lim_{x\to 0}\mathscr{F}_t \left[t^{2m}l_n(|t|,y)\right](x,y)=(-1)^{|m|}\lim_{x\to 0}\partial_x^{2m}\frac{\sinh(y\sqrt{|x|^2+\lambda^2})}{\sinh(a\sqrt{|x|^2+\lambda^2})},
\]
where
\[
\partial_x^{2m}=\frac{\partial^{2|m|}}{\partial x_1^{2m_1}\partial x_2^{2m_2}\dots\partial x_n^{2m_n}}.
\]
For $x\in\mathbb R^n$ we have the expansion
\[
L_n(x,y)=\frac{\sinh(y\sqrt{|x|^2+\lambda^2})}{\sinh(a\sqrt{|x|^2+\lambda^2})}=\sum_{|m|=0}^{\infty}p_{2|m|}(\lambda,y)\frac{(-1)^{|m|}}{|2m|!}(x_1^2+x_2^2\dots+x_n^2)^{|m|},
\]
the coefficient for $x^{2m}$ in this expansion is
\[
p_{2|m|}(\lambda,y)\frac{(-1)^{|m|}|m|!}{|2m|!m!},~~   m!=m_1!m_2!\dots m_n! , 
\]
Hence,
\[
p_{2m} (\lambda,y)=\frac{(2m)!|m|!}{|2m|!m!}p_{2|m|}(\lambda,y).
\]
For example,
\[
p_{2(2,1,1)}(\lambda,y)=\frac1{35}p_8 (\lambda,y).
\]
     The functions $v_k (x,y)= u_k (x,a-y)$ are solutions of the Helmholtz equation satisfying the boundary conditions
\[
v_k (x,0)=x^k,~~   v_k (x,a)=0,~~x\in\mathbb R^n,~~k~\text{is a multiindex}.
\]

\section{		Solution of the Dirichlet problem for the Helmholtz equation in the case $\nu=\mu^2$}

     Replacing $\lambda$ in (8), (9) by $i\mu,~\mu>0$, we obtain the Dirichlet problem
\[
\Delta u(x,y)+\mu^2 u(x,y)=0,~~\mu>0,~~x\in \mathbb R^n,~~0<y<a,   \eqno{(16)}
\]
\[
u(x,0)=\phi(x),~~u(x,a)=\psi(x),~~x\in\mathbb{R}^n, \eqno{(17)}
\]
where $\phi(x),\psi(x)$ are polynomials.

     All formulas obtained in section 3 are preserved with the replacement of hyperbolic functions by circular ones and the replacement of the minus sign by the plus sign on the right side of the recurrence formula (14). Namely, the solution of the Dirichlet problem (16), (17) is written in the form of a convolution
\[
                                       u(x,y)=l_n (|x|,a-y)*\phi(x)+l_n (|x|,y)*\psi(x),                           
\]
where
\[
  l_n (|x|,y)=\mathscr{F}_t^{-1}\left[L_n (|t|,y)\right](x,y),~~L_n(|t|,y)=\frac{\sin(y\sqrt{\mu^2-|t|^2})}{\sin(a\sqrt{\mu^2-|t|^2})}.
\]

If $\mu^2$ is not an eigenvalue of the Laplace operator in a layer with Dirichlet boundary conditions, then this solution will be unique in the class of functions of slow growth. 

For $\phi(x)=0,~~\psi(x)=x^k,~~x\in\mathbb R~~,k\in\mathbb N\cup\{0\}$, the solution to the Dirichlet problem is the function
\[
u_k(x,y)=\sum_{m=0}^{[k/2]}C_k^{2m}x^{k-2m}p_{2m}(\mu,y).
\]
where
\[
p_{2m}(\mu,y)=(2m-1)!!\left(\frac1{\mu}\frac{\partial}{\partial\mu}\right)^m\frac{\sin(\mu y)}{\sin(\mu a)},~~\mu\ne\pi k/a,~~k\in\mathbb N.
\]
For example, 
\[
u_2(x,y)=x^2\frac{\sin(\mu y)}{\sin(\mu a)}+\frac{y\cos(\mu y)}{\mu\sin(\mu a)}-\frac{a\sin(\mu y)\cos(\mu a)}{\mu\sin^2(\mu a)},~~\mu\ne\pi k/a,~~k\in\mathbb N.
\]
As $\mu$ tends to zero, $u_2 (x,y)$ go over into polynomial solution of the Dirichlet problem for the Laplace equation
\[
\lim_{\mu\to 0}u_2(x,y)=\frac{y}{3a}(3x^2-y^2+a^2).
\]

If $x\in\mathbb R^n,~~n>1,~~k~\text{is a multiindex}$, then the solutions are found by the same formulas in which now
\[
x^k=x_1^{k_1}x_2^{k_2}\dots x_n^{k_n},~~C_k^{2m}=C_{k_1}^{2m_1}C_{k_2}^{2m_2}\dots C_{k_n}^{2m_n},~~[k/2]=([k_1/2],[k_2/2],\dots,[k_n/2]),
\]
\[
p_{2m} (\mu,y)=\frac{(2m)!|m|!}{|2m|!m!}p_{2|m|}(\mu,y).
\]

\subsection{Uniqueness of the solution of the Dirichlet problem}

The solution of the Dirichlet problem (16), (17) will be unique in the class of functions of slow growth (the corresponding homogeneous boundary value problem has only trivial solutions in the class of functions of slow growth) if
\[
0<\mu<\frac{\pi}{a},~~0<\mu^2<\frac{\pi^2}{a^2}.
\]
If  $\mu^2\ge \pi^2/a,~\mu^2=\pi^2/a^2+b^2,~b^2\ge 0$, then the functions 
\[
\sin(b_1x_1)\sin(b_2x_2)\dots\sin(b_nx_n)\sin\left(\frac{\pi y}{a}\right),
\]
where $b_1,b_2,\dots,b_n$ are arbitrary non-negative numbers satisfying the equality
\[
 b_1^2+b_2^2+\dots+b_n^2=b^2,
\]
 are nontrivial solutions of slow growth  of the corresponding homogeneous boundary value problem.
 
We prove that in the case $0<\mu<\pi/a$ the corresponding homogeneous boundary value problem
\[
\Delta u(x,y)+\mu^2 u(x,y)=0,~~x\in \mathbb R^n,~~0<y<a,   
\]
\[
u(x,0)=0,~~u(x,a)=0,~~x\in\mathbb{R}^n, 
\]
has only trivial solutions in the class of functions of slow growth.

 The solution $u(x,y)$ as a function of the variable $y$ can be expanded in a Fourier series
\[
u(x,y)=\sum_{n=1}^{\infty}b_n(x)\sin\left(\frac{\pi ny}{a}\right),
\]
where the coefficients $b_n (x)$ are functions of slow growth. Substituting this function $u(x,y)$ into equation (16), we obtain the equations for the coefficients $b_n (x)$:
\[
\Delta b_n(x)+\left(\mu^2-\frac{\pi^2n^2}{a^2}\right)b_n(x)=0,~~n=1,2,\dots
\]

If $0<\mu<\pi/a,~0<\mu^2<\pi^2/a^2$, then $\mu^2-\pi^2n^2/a^2<0,~\forall n=1,2,\dots$ and these equations have only trivial solutions in the class of functions of slow growth.

\section{	Solution of the mixed Dirichlet-Neumann boundary value problem for the Helmoltz equation in the case $\nu=-\lambda^2$}

It is enough to consider a homogeneous equation.
\[
\Delta u(x,y)-\lambda^2 u(x,y)=0,~~\lambda>0,~~x\in \mathbb R^n,~~0<y<a,   \eqno{(19)}
\]
\[
u(x,0)=\phi(x),~~u_y(x,a)=\psi(x),~~x\in\mathbb{R}^n, \eqno{(20)}
\]
where $\phi(x),\psi(x)$ are polynomials.

                      As in the case of the Dirichlet problem, we have a unique solution in the class of functions of slow growth
\[
u(x,y)=l_n (|x|,y)*\phi(x)+k_n (|x|,y)*\psi(x),    
\]                      
where 
\[
  l_n (|x|,y)=\mathscr{F}_t^{-1}\left[L_n (|t|,y)\right](x,y),~~k_n(|x|,y) = \mathscr{F}_t^{-1}\left[K_n (|t|,y)\right](x,y),
\]
\[
L_n (|t|,y)=\frac{\cosh\left((a-y)\sqrt{|t|^2+\lambda^2}\right)}{\cosh\left(a\sqrt{|t|^2+\lambda^2}\right)},~~K_n (|t|,y)=\frac{\sinh\left(y\sqrt{|t|^2+\lambda^2}\right)}{\sqrt{|t|^2+\lambda^2}\cosh\left(a\sqrt{|t|^2+\lambda^2}\right)}.
\]

\subsection{	Case $n=1$}

If $\phi(x)=1,~\psi(x)=0,~x\in\mathbb R$, then the solution to problem (19), (20) is the function
\[
u_0(x,y)=\frac{\cosh\left((a-y)\lambda\right)}{\cosh(a\lambda)}.
\]
If $\phi(x)=0,~\psi(x)=1,~x\in\mathbb R$, then the solution to problem (19), (20) is the function
\[
v_0(x,y)=\frac{\sinh(y\lambda)}{\sinh(a\lambda)}.
\]
If $\phi(x)=x^k,~\psi(x)=0,~x\in\mathbb R,~k\in\mathbb N$, then the solution to problem (19), (20) is the function
\[
u_k(x,y)=\sum_{m=0}^{[k/2]}C_k^{2m}x^{k-2m}p_{2m}(\lambda,y),
\]
where the functions $p_{2m} (\lambda,y)$ are the coefficients of expansion in a power series in~$x$ of the function
\[
L_1 (x,y)=\frac{\cosh\left((a-y)\sqrt{x^2+\lambda^2}\right)}{\cosh\left(a\sqrt{x^2+\lambda^2}\right)}=\sum_{m=0}^{\infty}p_{2m}(\lambda,y)\frac{(-1)^mx^{2m}}{(2m)!},
\]
i.e. $L_1 (x,y)$ is the generating function for $p_{2m}(\lambda,y)$. These functions can be calculated using the recurrence formula
\[
p_0(\lambda y)=\frac{\cosh\left((a-y)\lambda\right)}{\cosh(a\lambda)},
\]
\[
p_{2m} (\lambda,y)=-(2m-1)\frac{1}{\lambda}\frac{\partial}{\partial \lambda}   p_{2m-2} (\lambda,y)),~~m=1,2,…      \eqno{(21)}
\] 
Hence,
\[
p_{2m}(\lambda,y)=(2m-1)!!\left(-\frac1{\lambda}\frac{\partial}{\partial\lambda}\right)^m\frac{\cosh\left((a-y)\lambda\right)}{\cosh(a\lambda)},
\]
and for example
\[
u_1(x,y)=x\cosh(\lambda y)-\frac{x\sinh(\lambda a)\sinh(\lambda y)}{\cosh(\lambda a)},
\]
\[
u_2(x,y)=x^2\cosh(\lambda y)-\frac{x^2\sinh(\lambda a)\sinh(\lambda y)}{\cosh(\lambda a)}+\frac{y\sinh(\lambda a)\cosh(\lambda y)}{\lambda\cosh(\lambda a)}+\frac{a\sinh(\lambda y)}{\lambda}-
\]
\[
-\frac{y\sinh(\lambda y)}{\lambda}-\frac{a\sinh^2(\lambda a)\sinh(\lambda y)}{\lambda\cos^2(\lambda a)}
\]
As $\lambda\to 0$, the functions $u_k (x,y)$ go over into polynomials that are solutions of the Dirichlet – Neumann boundary value problem for the Laplace equation \cite{Alg1}. For example,
\[
\lim_{\lambda\to 0}u_0(x,y)=1,~~\lim_{\lambda\to 0}u_1(x,y)=x,~~\lim_{\lambda\to 0}u_2(x,y)=x^2+2ay-y^2.
\]
If $\phi(x)=0,~\psi(x)=x^l,~x\in\mathbb R,~l\in\mathbb N$, then the solution to problem (19), (20) is the function
\[
v_l(x,y)=\sum_{m=0}^{[l/2]}C_l^{2m}x^{l-2m}q_{2m}(\lambda,y),
\]
where the functions $q_{2m}(\lambda,y)$ are the coefficients of expansion in a power series in~$x$ of the function
\[
K_1 (x,y)=\frac{\sinh\left(y\sqrt{x^2+\lambda^2}\right)}{\sqrt{x^2+\lambda^2}\cosh\left(a\sqrt{x^2+\lambda^2}\right)}=\sum_{m=0}^{\infty}q_{2m}(\lambda,y)\frac{(-1)^mx^{2m}}{(2m)!},
\]
i.e. $K_1(x,y)$ is the generating function for $q_{2m}(\lambda,y)$. These functions can be calculated using the recurrece formula
\[
q_0(\lambda,y)=\frac{\sinh(\lambda y)}{\lambda\cosh(\lambda a)},
\]
\[
q_{2m} (\lambda,y)=-(2m-1)\frac{1}{\lambda}\frac{\partial}{\partial \lambda}q_{2m-2}(\lambda,y)),~~m=1,2,…      \eqno{(22)}
\] 
Hence,
\[
q_{2m}(\lambda,y)=(2m-1)!!\left(-\frac1{\lambda}\frac{\partial}{\partial\lambda}\right)^m\frac{\sinh(y\lambda)}{\lambda\cosh(a\lambda)},
\]
and for example,
\[
v_1(\lambda,y)=\frac{x\sinh(\lambda y)}{\lambda\cosh(\lambda a)},
\]
\[
v_2(\lambda,y)=\frac{x^2\sinh(\lambda y)}{\lambda\cosh(\lambda a)}-\frac{y\cosh(\lambda y)}{\lambda^2\cosh(\lambda a)}+\frac{\sinh(\lambda y)}{\lambda^3\cosh(\lambda a)}+\frac{a\sinh(\lambda a)\sinh(\lambda y)}{\lambda^2\cosh^2(\lambda a)}.
\]
As $\lambda\to 0$, the functions $v_l (x,y)$ go over into polynomials that are solutions of the Dirichlet – Neumann boundary value problem for the Laplace equation \cite{Alg1}. For example,
\[
\lim_{\lambda\to 0}v_0(x,y)=y,~~\lim_{\lambda\to 0}v_1(x,y)=xy,~~\lim_{\lambda\to 0}v_2(x,y)=x^2y-\frac{y^3}{3}+a^2y.
\]
\textbf{Example 3. }
Consider the Dirichlet-Neumann problem for the inhomogeneous Helmholtz equation
\[
\Delta u(x,y)-\lambda^2 u(x,y)=x^2y^2,~~-\infty<x<\infty,~~0<y<a,~~\lambda>0,   
\]
\[
u(x,0)=0,~~u_y(x,a)=0,~~-\infty<x<\infty.
\]
Just as in the case of the Dirichlet problem, this problem reduces to the problem for a homogeneous equation
\[
\Delta v(x,y)-\lambda^2 v(x,y)=0,~~-\infty<x<\infty,~~0<y<a,~~\lambda>0,   
\]
\[
v(x,0)=\frac{2x^2}{\lambda^4}+\frac{8}{\lambda^6},~~v_y(x,a)=\frac{2x^2a}{\lambda^2}+\frac{4a}{\lambda^4},~~-\infty<x<\infty.
\]
The solution to the original problem is the function
\[
u(x,y)=-\frac{x^2y^2}{\lambda^2}-\frac{2y^2}{\lambda^4}-\frac{2x^2}{\lambda^4}-\frac{8}{\lambda^6}+\frac{2}{\lambda^4}u_2(x,y)+
\]
\[
+\frac{8}{\lambda^6}u_0(x,y)+\frac{2a}{\lambda^2}v_2(x,y)+\frac{4a}{\lambda^4}v_0(x,y)=-\frac{x^2y^2}{\lambda^2}-\frac{2y^2}{\lambda^4}-\frac{2x^2}{\lambda^4}-\frac{8}{\lambda^6}+
\]
\[
+\frac{2x^2\cosh(\lambda y)}{\lambda^4}-\frac{2x^2\sinh(\lambda a)\sinh(\lambda y)}{\lambda^4\cosh(\lambda a)}+\frac{2y\sinh(\lambda a)\cosh(\lambda y)}{\lambda^5\cosh(\lambda a)}+
\]
\[
+\frac{2a\sinh(\lambda y)}{\lambda^5}-\frac{2y\sinh(\lambda y)}{\lambda^5}-\frac{2a\sinh^2(\lambda a)\sinh(\lambda y)}{\lambda^5\cosh^2(\lambda a)}+\frac{8\cosh(\lambda (a-y))}{\lambda^6\cosh(\lambda a)}+
\]
\[
+\frac{2ax^2\sinh(\lambda y)}{\lambda^3\cosh(\lambda a)}-\frac{2ay\cosh(\lambda y)}{\lambda^4\cosh(\lambda a)}+\frac{6a\sinh(\lambda y)}{\lambda^5\cosh(\lambda a)}+\frac{2a^2\sinh(\lambda a)\sinh(\lambda y)}{\lambda^4\cosh^2(\lambda a)}.
\]
As $\lambda\to 0$, this function go over into a polynomial,
\[
\lim_{\lambda\to 0}u(x,y)=\frac{1}{12}x^2y^4-\frac{a^3}{3}x^2y-\frac{1}{180}y^6+\frac{a^3}{9}y^3-\frac{3a^5}{10}y,
\]
which is a solution of the Dirichlet-Neumann boundary value problem for the Poisson equation \cite{Alg1}
\[
\Delta u(x,y)=x^2y^2,~~-\infty<x<\infty,~~0<y<a,   
\]
\[
u(x,0)=0,~~u_y(x,a)=0,~~-\infty<x<\infty.
\]

\subsection{	Case $n>1$}

Similarly to the case of the Dirichlet problem, we obtain that the solution of the Dirichlet-Neumann problem for $\phi(x)=x^k,~\psi(x)=0,~x\in\mathbb R^n,~k=(k_1,\dots,k_n )$, is a function
\[
u_k(x,y)=\sum_{m=0}^{[k/2]}C_k^{2m}x^{k-2m}p_{2m}(\lambda,y),
\]
where
\[
x^k=x_1^{k_1}x_2^{k_2}\dots x_n^{k_n},~~C_k^{2m}=C_{k_1}^{2m_1}C_{k_2}^{2m_2}\dots C_{k_n}^{2m_n},~~[k/2]=([k_1/2],[k_2/2],\dots,[k_n/2]),
\]
\[
p_{2m} (\lambda,y)=\frac{(2m)!|m|!}{|2m|!m!}p_{2|m|}(\lambda,y).
\]
For $\phi(x)=0,~\psi(x)=x^k$, the solution to the Dirichlet-Neumann problem is the function
\[
v_k(x,y)=\sum_{m=0}^{[k/2]}C_k^{2m}x^{k-2m}q_{2m}(\lambda,y),~q_{2m} (\lambda,y)=\frac{(2m)!|m|!}{|2m|!m!}q_{2|m|}(\lambda,y).
\]
The functions $p_{2|m|}(\lambda,y)$ and $q_{2|m|}(\lambda,y)$ are found by formulas (21) and (22), respectively.

\section{	The solution of the mixed Dirichlet-Neumann boundary value problem for the Helmholtz equation in the case $\nu=~\mu^2$}

      All formulas obtained in section 5 are preserved with the replacement of hyperbolic functions by circular ones and the replacement of the minus sign by the plus sign on the right side of recurrence formulas (21), (22). Namely, the solution of the Dirichlet-Neumann problem
\[
\Delta u(x,y)+\mu^2 u(x,y)=0,~~\mu>0,~~x\in \mathbb R^n,~~0<y<a,   \eqno{(23)}
\]
\[
u(x,0)=\phi(x),~~u_y(x,a)=\psi(x),~~x\in\mathbb{R}^n, \eqno{(24)}
\]
is written as a convolution
\[
u(x,y)=l_n (|x|,y)*\phi(x)+k_n (|x|,y)*\psi(x),    
\]                      
where 
\[
  l_n (|x|,y)=\mathscr{F}_t^{-1}\left[L_n (|t|,y)\right](x,y),~~k_n(|x|,y) = \mathscr{F}_t^{-1}\left[K_n (|t|,y)\right](x,y),
\]
\[
L_n (|t|,y)=\frac{\cos\left((a-y)\sqrt{\mu^2-|t|^2}\right)}{\cos\left(a\sqrt{\mu^2-|t|^2}\right)},~~K_n (|t|,y)=\frac{\sin\left(y\sqrt{\mu^2-|t|^2}\right)}{\sqrt{\mu^2-|t|^2}\cos\left(a\sqrt{\mu^2-|t|^2}\right)}.
\]
      If $\mu^2$ is not an eigenvalue of the Laplace operator in a layer with Dirichlet-Neumann boundary conditions, then this solution will be unique in the class of functions of slow growth. 

 For $\phi(x)=x^k,~~\psi(x)=0,~~x\in\mathbb R~~,k\in\mathbb N\cup\{0\}$, the solution to the Dirichlet-Neumann problem is the function
\[
u_k(x,y)=\sum_{m=0}^{[k/2]}C_k^{2m}x^{k-2m}p_{2m}(\mu,y).
\]
where
\[
p_{2m}(\mu,y)=(2m-1)!!\left(\frac1{\mu}\frac{\partial}{\partial\mu}\right)^m\frac{\cos(\mu (a-y))}{\cos(\mu a)},~~\mu\ne\pi/2a+\pi k/a,~~k\in\mathbb N.
\]

For $\phi(x)=0,~~\psi(x)=x^k,~~x\in\mathbb R~~,k\in\mathbb N\cup\{0\}$, the solution to the Dirichlet-Neumann problem is the function
\[
v_k(x,y)=\sum_{m=0}^{[k/2]}C_k^{2m}x^{k-2m}q_{2m}(\mu,y).
\]
where
\[
q_{2m}(\mu,y)=(2m-1)!!\left(\frac1{\mu}\frac{\partial}{\partial\mu}\right)^m\frac{\sin(\mu y)}{\mu\cos(\mu a)},~~\mu\ne\pi/2a+\pi k/a,~~k\in\mathbb N.
\]
For example, 
\[
v_2(x,y)=x^2\frac{\sin(\mu y)}{\cos(\mu a)}-\frac{y\cos(\mu y)}{\mu^2\cos(\mu a)}+\frac{\sin(\mu y)}{\mu^3\cos(\mu a)}+\frac{a\sin(\mu a)\sin(\mu y)}{\mu^2\cos^2(\mu a)},~~\mu\ne\pi/2a+\pi k/a,~~k\in\mathbb N.
\]
As $\mu$ tends to zero, $v_2 (x,y)$ go over into polynomial solution of the Dirichlet-Neumann problem for the Laplace equation
\[
\lim_{\mu\to 0}v_2(x,y)=x^2y-\frac{y^3}{3}+a^2y.
\]

If $x\in\mathbb R^n,~~n>1,~~k,m~\text{are  multiindeces}$, then the solutions are found by the same formulas in which now
\[
x^k=x_1^{k_1}x_2^{k_2}\dots x_n^{k_n},~~C_k^{2m}=C_{k_1}^{2m_1}C_{k_2}^{2m_2}\dots C_{k_n}^{2m_n},~~[k/2]=([k_1/2],[k_2/2],\dots,[k_n/2]),
\]
\[
p_{2m} (\mu,y)=\frac{(2m)!|m|!}{|2m|!m!}p_{2|m|}(\mu,y),~q_{2m} (\mu,y)=\frac{(2m)!|m|!}{|2m|!m!}q_{2|m|}(\mu,y).
\]
 
\subsection{Uniqueness of the solution of the Dirichlet-Neumann problem}  

      Similarly to the case of the Dirichlet problem,  the solution of the Dirichlet-Neumann problem (23), (24) will be unique in the class of functions of slow growth if 
\[
0<\mu<\frac{\pi}{2a},~0<\mu^2<\frac{\pi^2}{4a^2}.
\]

If $\mu^2\ge\pi^2 /4a^2,~\mu^2=\pi^2 /4a^2+b^2,~b^2\ge 0$,  then functions 
\[
\sin(b_1x_1)\sin(b_2x_2)\dots\sin(b_nx_n)\sin\left(\frac{\pi y}{2a}\right),
\]
where $b_1,b_2,\dots,b_n$ are arbitrary non-negative numbers satisfying the equality
\[
b_1^2+b_2^2+\dots+b_n^2=b^2,
\]
are nontrivial solutions of slow growth of the corresponding homogeneous boundary value problem.

\textbf{Example 4.}
Consider the Dirichlet-Neumann problem for the inhomogeneous Helmholtz equation
\[
\Delta u(x,y)+\mu^2u(x,y)=x^{(2,1,1)}y^3,~~x\in\mathbb R^3,~~0<y<1,~~0<\mu<\pi/2,
\]
\[
u(x,0)=0,~~u(x,1)=0,~~x\in\mathbb R^3.
\]

By formula (7) we obtain a polynomial solution of the inhomogeneous equation
\[
\tilde u(x,y)=\frac{1}{\mu^2}x^{(2,1,1)}y^3-\frac{2}{\mu^4}x^{(0,1,1)}y^3-\frac{6}{\mu^4}x^{(2,1,1)}y+\frac{24}{\mu^6}x^{(0,1,1)}y,
\]
and the problem reduces to the problem for a homogeneous equation
\[
\Delta v(x,y)+\mu^2v(x,y)=0,~~x\in\mathbb R^3,~~0<y<1,~~0<\mu<\pi/2,
\]
\[
v(x,0)=-\tilde u(x,0)=0,~~~~x\in\mathbb R^3,
\]
\[
v_y(x,1)=-\tilde u_y(x,1)=-\frac{3}{\mu^2}x^{(2,1,1)}+\frac{6}{\mu^4}x^{(0,1,1)}+\frac{6}{\mu^4}x^{(2,1,1)}-\frac{24}{\mu^6}x^{(0,1,1)},~~x\in\mathbb R^3.
\]
The solution to this problem (which is unique in the class of functions of slow growth) is
\[
v(x,y)=-\frac{3}{\mu^2}v_{(2,1,1)}(x,y)+\frac{6}{\mu^4}v_{(0,1,1)}(x,y)+\frac{6}{\mu^4}v_{(2,1,1)}(x,y)-\frac{24}{\mu^6}v_{(0,1,1)}(x,y).
\]
The solution to the original problem is a quasipolynomial
\[
u(x,y)=\tilde u(x,y)+v(x,y)=
\]
\[
=\left(\frac{y^3}{\mu^2}-\frac{6y}{\mu^4}+\frac{6\sin(\mu y)}{\mu^5\cos(\mu)}-\frac{3\sin(\mu y)}{\mu^3\cos(\mu)}\right)x^{(2,1,1)}+
\]
\[
+\left(-\frac{2y^3}{\mu^4}+\frac{24y}{\mu^6}+\frac{9\sin(\mu y)}{\mu^5\cos(\mu)}-\frac{30\sin(\mu y)}{\mu^7\cos(\mu)}+\frac{6y\cos(\mu y)}{\mu^6\cos(\mu)}\right)x^{(0,1,1)}+
\]
\[
+\left(-\frac{3y\cos(\mu y)}{\mu^4\cos(\mu)}+\frac{6\sin(\mu y)\sin(\mu)}{\mu^6\cos^2(\mu)}-\frac{3\sin(\mu y)\sin(\mu)}{\mu^4\cos^2(\mu)}\right) x^{(0,1,1)}.
\]
As $\mu\to 0$, this solution go over into a polynomial
\[
\lim_{\mu\to 0}u(x,y)=\left(\frac{1}{20}y^5-\frac1{4}y\right)x^{(2,1,1)}+\left(-\frac{1}{420}y^7+\frac{1}{12}y^3-\frac{7}{30}y\right)x^{(0,1,1)},
\]
which is the solution of the Dirichlet-Neumann problem for the Poisson equation
\[
\Delta u(x,y)=x^{(2,1,1)}y^3,~~x\in\mathbb R^3,~~0<y<1,
\]
\[
u(x,0)=0,~~u(x,1)=0,~~x\in\mathbb R^3.
\]

\section*{Conclusion}

The paper considers the Dirichlet and Dirichlet-Neumann problems in a multidimensional infinite layer for the inhomogeneous Helmholtz equation with polynomial right-hand side and polynomials in the right-hand sides of the boundary conditions. It is shown that the unique solution to each of these problems in the class of functions of slow growth in the case when the equation parameter is not an eigenvalue is a quasipolynomial. An algorithm for constructing this quasipolynomial is given and examples are considered.

\end{document}